# THE EMPIRICAL PROCESS ON GAUSSIAN SPHERICAL HARMONICS

By Domenico Marinucci[1] and Mauro Piccioni

*Università di Roma "Tor Vergata" and Università di Roma "La Sapienza"*

We establish weak convergence of the empirical process on the spherical harmonics of a Gaussian random field in the presence of an unknown angular power spectrum. This result suggests various Gaussianity tests with an asymptotic justification. The issue of testing for Gaussianity on isotropic spherical random fields has recently received strong empirical attention in the cosmological literature, in connection with the statistical analysis of cosmic microwave background radiation.

**1. Introduction.** In recent years an enormous amount of attention has been devoted to testing for Gaussianity for spherical random fields, especially in the astrophysical and cosmological literature. The empirical motivation for these studies can be explained as follows. The ongoing NASA satellite mission MAP and the forthcoming ESA mission Planck will probe cosmic microwave background radiation (CMB) to an unprecedented accuracy. CMB can be viewed as a signature of the distribution of matter and radiation in the very early universe, and as such it is expected to yield very tight constraints on physical models for the Big Bang. For the density fluctuations of this field, the highly popular inflationary scenario [see Peebles (1993) or Peacock (1999)] predicts a Gaussian distribution, whereas alternative cosmological theories, such as topological defects or nonstandard inflationary models, predict other types of behavior. The density distributions of fluctuations are also instrumental for drawing correct inferences on the physical constants which can be estimated from CMB radiation; indeed, point and interval estimation procedures for cosmological parameters have been based almost exclusively upon Gaussian assumptions, which of course

Received July 2002; revised April 2003.
[1]Supported by MURST.

need to be validated before reliable statistical inference can take place. For these reasons, many different Gaussianity tests were considered in the recent cosmological literature, some of them based upon the topological properties of Gaussian fields [Novikov, Schmalzing and Mukhanov (2000), Phillips and Kogut (2001) and Dorè, Colombi and Bouchet (2003)], others on higher-order cumulant spectra [Winitzki and Wu (2000) and Komatsu and Spergel (2001)].

Let $T(\theta, \varphi)$ be a random field indexed by the unit sphere $S^2$; that is, for each azimuth $0 \leq \theta \leq \pi$ and elongation $0 \leq \varphi < 2\pi$, $T(\theta, \varphi)$ is a random variable defined on some probability space. Throughout this paper we assume that $T(\theta, \varphi)$ has zero mean, finite variance, is mean square continuous, and is isotropic, that is, its covariance is invariant with respect to the group of rotations. Furthermore, we introduce the spherical harmonics, defined by

$$Y_{l,m}(\theta, \varphi) = \begin{cases} \sqrt{\dfrac{2l+1}{4\pi} \dfrac{(l-m)!}{(l+m)!}} P_{lm}(\cos\theta) \exp(im\varphi), & \text{for } m > 0, \\ (-1)^m Y^*_{l,-m}(\theta, \varphi), & \text{for } m < 0, \end{cases}$$

where the asterisk denotes complex conjugation and $P_{lm}(\cos\theta)$ denotes the associated Legendre polynomial of degree $l$, $m$, that is,

$$P_{lm}(x) = (-1)^m (1-x^2)^{m/2} \frac{d^m}{dx^m} P_l(x),$$

$$P_l(x) = \frac{1}{2^l l!} \frac{d^l}{dx^l} (x^2 - 1)^l, \qquad m = 0, 1, 2, \ldots, l, l = 1, 2, 3, \ldots.$$

A detailed discussion of the properties of the spherical harmonics can be found, for instance, in Liboff (1998), Chapter 9, and in Varshalovich, Moskalev and Khersonskii (1988), Chapter 5. The following spectral representation holds in the mean square sense [see, e.g., Hannan (1970), Wong (1971) and Leonenko (1999)]:

$$(1) \qquad T(\theta, \varphi) = \sum_{l=1}^{\infty} \sum_{m=-l}^{l} a_{lm} Y_{lm}(\theta, \varphi),$$

where the triangular array $\{a_{lm}\}$ represents a set of random coefficients which can be obtained from $T(\theta, \varphi)$ through the inversion formula

$$(2) \qquad a_{lm} = \int_0^{2\pi} \int_0^{\pi} T(\theta, \varphi) Y^*_{lm}(\theta, \varphi) \sin\theta \, d\theta \, d\varphi,$$

$$m = 0, \pm 1, \ldots, \pm l, l = 1, 2, \ldots.$$

These coefficients are zero-mean and uncorrelated [see Wong (1971), pages 253 and 254]; hence, if $T(\theta, \varphi)$ is Gaussian, as we shall always assume in this paper, they have a complex Gaussian distribution, and they are independent over $l$ and $m \geq 0$ (although $a_{l,-m} = (-1)^m a_{lm}^*$), with variance $E|a_{lm}|^2 = C_l$, $m = 0, \pm 1, \ldots, \pm l$. The sequence $\{C_l\}$ denotes the angular power spectrum of the random field; we shall always assume that $C_l$ is strictly positive for all values of $l$.

Our purpose in this paper is to construct Gaussianity tests based on the empirical distribution function for the triangular array of random coefficients $\{a_{lm}\}$, $l = 1, \ldots, L$, $m = 0, \pm 1, \ldots, \pm l$, as $L \to \infty$. In principle, expression (2) can be used to recover the coefficients $a_{lm}$ for any value of $l$. In practical applications, however, $L$ is finite and determined by the resolution of the experiment. For instance, the coefficients $a_{lm}$ may be contaminated by noise, with a decreasing signal to noise ratio as $l$ grows. Throughout this paper we adopt the simplifying assumption that the noise is negligible only for $l = 1, \ldots, L$, hence the higher-order coefficients are discarded. The current status of technology in satellite missions suggests that our assumption may provide a reasonable approximation for $L$ as large as several hundreds for the ongoing NASA experiment MAP, a value which is likely to increase to a few thousands for the forthcoming ESA experiment Planck.

Testing Gaussianity in a spherical field poses considerable extra difficulty vs. testing in Euclidean spaces, since (fixed radius) spherical fields are never ergodic. Our proposal is to consider the asymptotic behavior of the empirical process on the triangular array $\{a_{lm}\}$, in the presence of infinitely many unknown parameters $\{C_l\}$; in fact, the angular power spectrum is typically unknown in practice, and needs to be estimated nonparametrically from the data [Miller, Nichol, Genovese and Wasserman (2002) and Wasserman, Miller, Nichol, Genovese, Jang, Connolly, Moore, Schneider and the PICA Group (2001)]. The presence of a growing sequence of estimated parameters makes it hard to exploit the modern theory of empirical processes, as presented for instance by van der Vaart and Wellner (1996) or Dudley (1999); we hence resort to a more traditional approach, based upon standard weak convergence theorems in the Skorohod space $D[0,1]^2$, as presented for instance by Bickel and Wichura (1971) or Shorack and Wellner (1986). We refer to these works for the definition and topological properties of such a space. For the sequel it is enough to recall that, in order to prove the weak convergence of a sequence $\{K_L(r, \alpha)\}$ of $D[0,1]^2$-valued processes to the field $K(r, \alpha)$, we need to prove both convergence of all finite-dimensional distributions and tightness. For the latter we will use a sufficient criterion due to Bickel and Wichura (1971).

First, for a generic "block" $B = (\alpha_1, \alpha_2] \times (r_1, r_2] \subset [0,1]^2$, define the increments of the field $K_L$,

$$K_L(B) = K_L(\alpha_2, r_2) - K_L(\alpha_2, r_1) - K_L(\alpha_1, r_2) + K_L(\alpha_1, r_1).$$

We define two types of adjacent blocks:

Type I blocks

$$B_1 = (\alpha_1, \alpha] \times (r_1, r_2],$$
$$B_2 = (\alpha, \alpha_2] \times (r_1, r_2], \qquad 0 \leq \alpha_1 \leq \alpha \leq \alpha_2 \leq 1, 0 \leq r_1 \leq r_2 \leq 1,$$

and Type II blocks

$$B_1 = (\alpha_1, \alpha_2] \times (r_1, r],$$
$$B_2 = (\alpha_1, \alpha_2] \times (r, r_2], \qquad 0 \leq \alpha_1 \leq \alpha_2 \leq 1, 0 \leq r_1 \leq r \leq r_2 \leq 1.$$

Bickel and Wichura (1971, Theorem 3) show that tightness of the sequence $\{K_L\}$ is satisfied if there exist $\beta > 1$, $\gamma > 0$ such that for all blocks $B_1$ and $B_2$,

$$E([\min\{|K_L(B_1)|, |K_L(B_2)|\}]^\gamma) \leq C(\mu(B_1 \cup B_2))^\beta,$$

which is implied by the stronger condition

(3) $\quad E(|K_L(B_1)|^{\gamma_1} |K_L(B_2)|^{\gamma_2}) \leq C(\mu(B_1 \cup B_2))^\beta, \qquad \gamma_1 + \gamma_2 = \gamma,$

where $\mu$ is some finite measure on $[0,1]^2$ with continuous marginals. They also show that for some particular class of processes (partial-sum processes) we can restrict to blocks having corners in $[0,1] \times \{0, \frac{1}{L}, \ldots, 1\}$, so that $(r_2 - r), (r_1 - r) \geq \frac{1}{L}$ always [see Bickel and Wichura (1971), page 1665].

The plan of this paper is as follows. In Section 2 we start from the case where the angular spectrum is known and then provide a formal definition of the empirical process for spherical harmonics we are interested in, for which we state the main weak convergence result. Section 3 presents the main steps of the proof; Section 4 draws some conclusions, which are also illustrated by a small Monte Carlo experiment, and points out some directions for further research. Many technical lemmas are collected separately in the Appendix. In the sequel, we use $C$ to denote a generic, positive and finite constant, whose value may vary from line to line. Also, throughout the paper we define $-\log x = \infty$ and (by continuity) $x \log x = x \log^2 x = 0$ for $x = 0$.

**2. The empirical process with unknown angular power spectrum.** We start by assuming that the sequence of coefficients $\{C_l\}_{l=1,2,\ldots}$ in the angular power spectrum of $T(\theta, \varphi)$ is known. Now recall that, under Gaussianity, $|a_{l0}|^2/C_l$ and $\{2|a_{lm}|^2\}/C_l = \{2|a_{l,-m}|^2\}/C_l$ are mutually independent chi-square variables, with one and two degrees of freedom, respectively. The special distributional properties of the single term $|a_{l0}|^2$ greatly complicates notation, whereas it can be shown that this term has no effect on asymptotic distributions; hence in the sequel we shall simply drop it and focus on $|a_{lm}|^2$ for $m = 1, 2, \ldots, l$. First introduce the Smirnov transformation

(4) $\quad u_{lm} = 1 - \exp\left(-\frac{|a_{lm}|^2}{C_l}\right), \qquad m = 1, 2, \ldots, l, l = 1, 2, \ldots, L,$

to convert the random variables $a_{lm}$ to a triangular array $\{u_{lm}\}$ of i.i.d. random variables with a uniform distribution in $[0,1]$. We can hence define their empirical distribution function over the $l$th row,

$$F_l(\alpha) = \frac{1}{l}\sum_{m=1}^{l} \mathbb{1}(u_{lm} \leq \alpha), \qquad 0 \leq \alpha \leq 1$$

[$\mathbb{1}(\cdot)$ denoting the indicator function], and the empirical process

$$G_l(\alpha) = \sqrt{l}\{F_l(\alpha) - \alpha\} = \frac{1}{\sqrt{l}}\left\{\sum_{m=1}^{l}\mathbb{1}(u_{lm} \leq \alpha) - \alpha\right\}, \qquad 0 \leq \alpha \leq 1.$$

In order to detect departures from the Gaussianity assumptions over some region of the angular decomposition, we shall consider the doubly indexed, integrated empirical process

$$K_L(r,\alpha) \stackrel{\text{def}}{=} \frac{1}{\sqrt{L}}\sum_{l=1}^{[Lr]} G_l(\alpha).$$

$G_l(\alpha)$ is zero-mean with covariance function $(\alpha_1 \wedge \alpha_2)(1 - (\alpha_1 \vee \alpha_2))$, as is well known. The integrated process $K_L(r,\alpha)$ has independent increments in $r$, from which it is easily obtained that its limiting covariance function should be

$$\lim_{L\to\infty} E K_L(r_1,\alpha_1) K_L(r_2,\alpha_2) = (r_1 \wedge r_2)(\alpha_1 \wedge \alpha_2)(1 - (\alpha_1 \vee \alpha_2)).$$

The Gaussian zero-mean process $K$ with this covariance function is called the Kiefer–Müller process on $[0,1]^2$. In fact it is a standard result that $K_L$, whose sample paths clearly belong to $D[0,1]^2$, converges weakly to this field as $L \to \infty$. We will be concerned, however, with the much more interesting (and difficult) case when the power spectrum is unknown and we have to estimate it nonparametrically from the data. A natural candidate to replace $C_l$ is the maximum likelihood estimate

$$\widehat{C}_l = \frac{1}{l}\sum_{m=1}^{l} |a_{lm}|^2.$$

The $u_{lm}$ are then replaced by the random variables

$$y_{lm} = 1 - \exp\left(-\frac{|a_{lm}|^2}{\widehat{C}_l}\right) = 1 - \exp(-l\xi_{lm}),$$

$$\xi_{lm} = \frac{|a_{lm}|^2}{\sum_{k=1}^{l}|a_{lk}|^2}, \qquad m = 1, 2, \ldots, l,$$

where $(\xi_{l1},\ldots,\xi_{ll})$ no longer has independent components, but has a Dirichlet distribution with parameters $(1,\ldots,1)$, that is, it is uniformly distributed on the unit simplex of $R^l$. For the sequel we recall that the $\xi_{li}$, $i = 1, \ldots, l$, are exchangeable with marginal distributions given by

$$P\{\xi_{l1} \leq \alpha\} = 1 - (1-\alpha)^{l-1}, \qquad \alpha \in [0,1], \tag{5}$$

and

$$P\{\xi_{l1} \leq \alpha_1, \xi_{l2} \leq \alpha_2\} = 1 - (1-\alpha_1)^{l-1} - (1-\alpha_2)^{l-1} + (1-\alpha_1-\alpha_2)^{l-1}, \tag{6}$$

$$\alpha_1, \alpha_2 \in [0,1].$$

We define now

$$\widehat{G}_l(\alpha) = \frac{1}{\sqrt{l}} \left[ \sum_{m=1}^{l} \{\mathbb{1}(y_{lm} \leq \alpha) - \alpha\} \right]$$

$$= \frac{1}{\sqrt{l}} \left[ \sum_{m=1}^{l} \left\{ \mathbb{1}\left( \xi_{lm} \leq -\frac{\log(1-\alpha)}{l} \right) - \alpha \right\} \right], \tag{7}$$

a representation which will be useful for the arguments to follow. Notice that, as usual, $\widehat{G}_l(0) = \widehat{G}_l(1) = 0$. The process (7) is known in the literature as the normalized uniform spacings process; see, for instance, Shorack and Wellner (1986), pages 731–733, where its limiting behavior as $l \to \infty$ is derived. More precisely, it is shown that $\widehat{G}_l$ converges weakly in $D[0,1]$ to $\widehat{G}_\infty$, where $\widehat{G}_\infty$ is a mean-zero Gaussian process with covariance function

$$E\widehat{G}_\infty(\alpha_1)\widehat{G}_\infty(\alpha_2) = (\alpha_1 \wedge \alpha_2)\{1 - (\alpha_1 \vee \alpha_2)\}$$
$$- (1-\alpha_1)(1-\alpha_2)\log(1-\alpha_1)\log(1-\alpha_2). \tag{8}$$

We shall focus instead on the partial sum of empirical processes

$$\widehat{K}_L(\alpha, r) = \frac{1}{\sqrt{L}} \sum_{l=1}^{[Lr]} \widehat{G}_l(\alpha), \qquad 0 \leq \alpha \leq 1, 0 \leq r \leq 1.$$

As we shall see, the asymptotic behavior of $\widehat{K}_L$ as $L \to \infty$ does not follow trivially from the asymptotics for $\widehat{G}_l$; intuitively, this is due to the effect of higher-order terms, which have to be controlled as they are summed over $l$. In particular we shall now show that the process $\widehat{K}_L$ has a nonnull asymptotic bias. Define

$$b_l(\alpha) = lE\left\{ \mathbb{1}\left( \xi_{l1} \leq -\frac{\log(1-\alpha)}{l} \right) - \alpha \right\}$$

so that

$$E\widehat{G}_l(\alpha) = \frac{1}{\sqrt{l}} b_l(\alpha).$$

LEMMA 2.1. *As $l \to \infty$,*

(9) $$\lim_{l \to \infty} b_l(\alpha) = (1-\alpha)\log(1-\alpha) + \tfrac{1}{2}(1-\alpha)\log^2(1-\alpha) = b(\alpha),$$

*and also, as $L \to \infty$,*

$$\lim_{L \to \infty} E\widehat{K}_L(\alpha, r) = 2\sqrt{r}\,b(\alpha).$$

PROOF. From (5) and for $l$ suitably large,

$$E\mathbb{1}\left(\xi_{l1} \leq \frac{-\log(1-\alpha)}{l}\right) = 1 - \left(1 + \frac{\log(1-\alpha)}{l}\right)^{l-1}, \qquad \alpha \in [0,1].$$

Therefore

$$b_l(\alpha) = l\left\{1 - \left(1 + \frac{\log(1-\alpha)}{l}\right)^{l-1} - \alpha\right\}.$$

Using Lemma A.2, we have immediately (9). Also

$$E\widehat{K}_L(x, r) = \frac{1}{\sqrt{L}} \sum_{l=1}^{[Lr]} \frac{1}{\sqrt{l}} \{b_l(\alpha) - b(\alpha)\} + \frac{b(\alpha)}{\sqrt{L}} \sum_{l=1}^{[Lr]} \frac{1}{\sqrt{l}}.$$

Now

$$\lim_{l \to \infty} \frac{b(\alpha)}{\sqrt{L}} \sum_{l=1}^{[Lr]} \frac{1}{\sqrt{l}} = b(\alpha)\sqrt{r}\int_0^1 \frac{1}{\sqrt{u}}\,du = 2b(\alpha)\sqrt{r}.$$

Also, because $\lim_{l \to \infty} b_l(\alpha) = b(\alpha)$, for any $\delta > 0$ there exists $l_0$ such that, for all $l > l_0$, $|b_l(\alpha) - b(\alpha)| < \delta$; hence

$$\frac{1}{\sqrt{L}} \sum_{l=1}^{[Lr]} \frac{1}{\sqrt{l}} \{b_l(\alpha) - b(\alpha)\} \leq \left\{\sup_{1 \leq l \leq l_0} |b_l(\alpha)| + |b(\alpha)|\right\} \frac{2\sqrt{l_0}}{\sqrt{L}} + 2\delta\sqrt{r} \leq \varepsilon,$$

for any $\varepsilon > 0$, for $L$ large enough, because $\delta$ is arbitrary. Thus the proof is completed. □

We shall hereafter write for brevity

$$\alpha_l = -\frac{\log(1-\alpha)}{l} \qquad \text{for } \alpha \in [0,1].$$

The limiting behavior of the covariances was anticipated in (8); indeed, we have the following.

LEMMA 2.2.   *For all $0 \leq \alpha_1,\ \alpha_2 \leq 1$,*

$$\lim_{l \to \infty} \text{Cov}\{\widehat{G}_l(\alpha_1), \widehat{G}_l(\alpha_2)\}$$

(10)
$$= (\alpha_1 \wedge \alpha_2)\{1 - (\alpha_1 \vee \alpha_2)\}$$
$$- (1 - \alpha_1)(1 - \alpha_2)\log(1 - \alpha_1)\log(1 - \alpha_2),$$

*and as $L \to \infty$, for all $0 \leq \alpha_1,\ \alpha_2 \leq 1,\ 0 \leq r_1,\ r_2 \leq 1$,*

$$\lim_{L \to \infty} \text{Cov}\{\widehat{K}_L(\alpha_1, r_1), \widehat{K}_L(\alpha_2, r_2)\}$$

(11)
$$= (r_1 \wedge r_2)[(\alpha_1 \wedge \alpha_2)\{1 - (\alpha_1 \vee \alpha_2)\}$$
$$- (1 - \alpha_1)(1 - \alpha_2)\log(1 - \alpha_1)\log(1 - \alpha_2)].$$

PROOF.   The first limiting result is standard in the theory of spacings [see Shorack and Wellner (1986) and Shorack (1972)]; indeed, from (5), (6) and Lemma A.1 the reader can check that

$$\lim_{l \to \infty} \text{Cov}\{\mathbb{1}(\xi_{l1} \leq \alpha_{1l}), \mathbb{1}(\xi_{l1} \leq \alpha_{2l})\}$$
$$= (\alpha_1 \wedge \alpha_2)\{1 - (\alpha_1 \vee \alpha_2)\},$$
$$\lim_{l \to \infty} l\,\text{Cov}\{\mathbb{1}(\xi_{l1} \leq \alpha_{1l}), \mathbb{1}(\xi_{l2} \leq \alpha_{2l})\}$$
$$= -\log(1 - \alpha_1)\log(1 - \alpha_2)(1 - \alpha_1)(1 - \alpha_2),$$

whence (10) follows easily. On the other hand, (11) is an immediate consequence of the Kronecker lemma and (10).   □

Now write

$$K_L^*(\alpha, r) = \widehat{K}_L(\alpha, r) - 2\sqrt{r}b(\alpha),$$

and define $K^*(\alpha, r)$ as the zero-mean Gaussian process on $[0,1] \times [0,1]$ with covariance

$$EK^*(\alpha_1, r_1)K^*(\alpha_2, r_2)$$
$$= (r_1 \wedge r_2)[(\alpha_1 \wedge \alpha_2)\{1 - (\alpha_1 \vee \alpha_2)\}$$
$$- (1 - \alpha_1)(1 - \alpha_2)\log(1 - \alpha_1)\log(1 - \alpha_2)].$$

To the best of our knowledge, the field $K^*(\alpha, r)$ has not appeared in the literature so far, and we label it a modified Kiefer–Müller process. Our main result will be the following.

THEOREM 2.1.   *As $L \to \infty$, weakly in $D[0,1]^2$,*

$$K_L^* \Rightarrow K^*,$$

$\Rightarrow$ *denoting weak convergence in the Skorohod space $D[0,1]^2$.*

The proof of this result will be given in the next section.

**3. The weak convergence proof.**

PROPOSITION 3.1. *As $L \to \infty$,*

$$K_L^* \stackrel{f.d.d.}{\longrightarrow} K^*.$$

PROOF. Consider $0 \le \alpha_1 < \cdots < \alpha_s \le 1$, $0 \le r_1 < \cdots < r_t \le 1$, and let

(12) $$G_l^*(\alpha, \beta] = \widehat{G}_l(\alpha, \beta] - E\widehat{G}_l(\alpha, \beta],$$

for $0 \le \alpha \le \beta \le 1$. Because of Lemma 2.1 the convergence of the $s \times t$ vector with components $K_L^*(\alpha_i, r_j)$ with $i = 1, 2, \ldots, s$, $j = 1, 2, \ldots, t$, to the same components of $K^*(\alpha, r)$ is equivalent to the joint convergence of the centered increments

$$\widehat{K}_L((\alpha_i, \alpha_{i+1}] \times (r_j, r_{j+1}]) - E\widehat{K}_L((\alpha_i, \alpha_{i+1}] \times (r_j, r_{j+1}])$$

$$= \frac{1}{\sqrt{L}} \sum_{l=[Lr_i]+1}^{[Lr_{i+1}]} [\widehat{G}_l(\alpha_{i+1}) - \widehat{G}_l(\alpha_i) - E\{\widehat{G}_l(\alpha_{i+1}) - \widehat{G}_l(\alpha_i)\}]$$

$$= \frac{1}{\sqrt{L}} \sum_{l=[Lr_i]+1}^{[Lr_{i+1}]} G_l^*(\alpha_i, \alpha_{i+1}]$$

for $i = 0, \ldots, s$ and $j = 0, \ldots, t$, where we have set $\alpha_0 = r_0 = 0$, $\alpha_{s+1} = r_{t+1} = 1$, to the increments

$$K^*((\alpha_i, \alpha_{i+1}] \times (r_j, r_{j+1}])$$
$$= K^*(\alpha_{i+1}, r_{j+1}) + K^*(\alpha_i, r_j) - K^*(\alpha_i, r_{j+1}) - K^*(\alpha_{i+1}, r_j).$$

Because of the independence over $l$ we can restrict ourselves to a fixed interval for $r$, by simplicity $(0, 1]$, say. It is now clear that we have to use the multidimensional central limit theorem (CLT) with independent but not identically distributed summands. A suitable control on the fourth moments clearly implies the Lindeberg condition, namely

(13) $$\frac{1}{L^2} \sum_{l=1}^{L} EG_l^*(\alpha_i, \alpha_{i+1}]^2 G_l^*(\alpha_k, \alpha_{k+1}]^2 = o(1) \quad \text{as } L \to \infty,$$

for any $i, k = 0, \ldots, t$, possibly equal. Now in Proposition 3.2 it will be shown that $EG_l^*(\alpha_i, \alpha_{i+1}]^2 G_l^*(\alpha_k, \alpha_{k+1}]^2$ is uniformly bounded by a constant, which will allow us to complete the proof. $\square$

PROPOSITION 3.2. *For $L = 1, 2, \ldots$, the sequence of fields $K_L^*(\alpha, r)$ is tight in $D[0, 1]^2$.*

Proof. We write
$$K_L^*(\alpha,r) = \widetilde{K}_L(\alpha,r) - E\widetilde{K}_L(\alpha,r) + \widehat{K}_L(\alpha,r)$$
$$- \widetilde{K}_L(\alpha,r) + E\widetilde{K}_L(\alpha,r) - 2\sqrt{r}b(\alpha),$$
where

(14) $$\widetilde{K}_L(\alpha,r) := \frac{1}{\sqrt{L}}\sum_{l=1}^{[Lr]}\widehat{G}_l(t_l(\alpha)),$$

$$t_l(\alpha) = \begin{cases} \alpha, & \text{for } \alpha < 1 - l^{-3/2}, \\ 1, & \text{for } \alpha \geq 1 - l^{-3/2}. \end{cases}$$

Notice that the modification $\widehat{G}_l(t_l(\alpha))$ differs from the integrated empirical process $\widehat{G}_l(\alpha)$ because it is tied down to zero at $\alpha = 1 - l^{-3/2}$, $l = 1, 2, \ldots, L$. The result will be clearly established if we can prove that the sequence $\widetilde{K}_L(\alpha,r) - E\widetilde{K}_L(\alpha,r)$ is tight, and as $L \to \infty$,

(15) $$\sup_{\alpha,r}|\widehat{K}_L(\alpha,r) - \widetilde{K}_L(\alpha,r)| = o_p(1),$$

(16) $$\sup_{\alpha,r}|E\widehat{K}_L(\alpha,r) - 2\sqrt{r}b(\alpha)| = o(1).$$

The proofs of (15) and (16) are given in Lemmas A.4 and A.5, respectively. We shall hence focus on tightness; let us write
$$\widetilde{K}'_L(\alpha,r) = \widetilde{K}_L(\alpha,r) - E\widetilde{K}_L(\alpha,r).$$
It is sufficient to establish that, for some probability measure $\mu(\cdot)$ with continuous marginals,
$$E(|\widetilde{K}'_L(B_1)|^2|\widetilde{K}'_L(B_2)|^2) \leq C(\mu(B_1 \cup B_2))^2,$$
where $B_1, B_2$ are either Type I or Type II blocks. Let us consider Type II blocks first, for which

(17) $$E\{\widetilde{K}'_L((\alpha_1,\alpha_2] \times (r_1,r])^2\widetilde{K}'_L((\alpha_1,\alpha_2] \times (r,r_2])^2\}$$
$$= E\{\widetilde{K}'_L((\alpha_1,\alpha_2] \times (r_1,r])^2\}E\{\widetilde{K}'_L((\alpha_1,\alpha_2] \times (r,r_2])^2\}.$$

Now recalling (12),
$$E\{\widetilde{K}'_L((\alpha_1,\alpha_2] \times (r_1,r])^2\} = \frac{1}{L}\sum_{l=[Lr_1]+1}^{[Lr]} EG_l^{*2}((\alpha_1,\alpha_2] \cap (0, 1-l^{-3/2}]).$$

For $0 \leq \alpha \leq \beta \leq 1$,
$$G_l^*(\alpha,\beta] = \frac{1}{\sqrt{l}}\sum_{m=1}^{l} Z_{lm}(\alpha,\beta],$$

where $Z_{lm}(\alpha, \beta] = Z_{lm}(\beta) - Z_{lm}(\alpha)$ and

$$Z_{lm}(\alpha) = \mathbb{1}\left(\xi_{lm} \leq -\frac{\log(1-\alpha)}{l}\right) - E\mathbb{1}\left(\xi_{lm} \leq -\frac{\log(1-\alpha)}{l}\right).$$

Thus

$$EG_l^{*2}((\alpha_1, \alpha_2] \cap (0, 1 - l^{-3/2}])$$
$$= \begin{cases} 0, & \text{for } \alpha_1 > 1 - l^{-3/2}, \\ EG_l^{*2}((\alpha_1, \alpha_2 \wedge (1 - l^{-3/2})]), & \text{otherwise.} \end{cases}$$

For $0 \leq \alpha \leq 1$ write $\tau_l(\alpha) = \alpha \wedge (1 - l^{-3/2})$, so that

$$EG_l^{*2}(\alpha_1, \tau_l(\alpha_2)]$$
$$= \frac{1}{l} \sum_{m_1=1}^{l} \sum_{m_2=1}^{l} EZ_{lm_1}(\alpha_1, \tau_l(\alpha_2)] Z_{lm_2}(\alpha_1, \tau_l(\alpha_2)]$$
$$= EZ_{l1}^2(\alpha_1, \tau_l(\alpha_2)] + (l-1) EZ_{l1}(\alpha_1, \tau_l(\alpha_2)] Z_{l2}(\alpha_1, \tau_l(\alpha_2)].$$

Now note that

$$EZ_{l1}^2(\alpha, \tau_l(\alpha_2)]$$
$$= \text{Var}(\mathbb{1}(\alpha_1 < y_{lm} \leq \tau_l(\alpha_2))) \leq E\mathbb{1}(\alpha_1 < y_{lm} \leq \tau_l(\alpha_2))$$
$$= p_l(\alpha_1, \tau_l(\alpha_2)) \leq e^2 |\tau_l(\alpha_2) - \alpha_1| \leq C|\alpha_2 - \alpha_1|,$$

the last step following from Lemma A.3. Also, from Lemma A.6,

$$(l-1)EZ_{l1}(\alpha_1, \tau_l(\alpha_2)] Z_{l2}(\alpha_1, \tau_l(\alpha_2)] \leq Cq(\alpha_1, \tau_l(\alpha_2)) \leq Cq(\alpha_1, \alpha_2)$$

with

$$q(\alpha_1, \alpha_2) = \int_{\alpha_1}^{\alpha_2} (1 + |\log(1-y)|) \, dy;$$

indeed,

$$q^2(a, b) \leq Cq(a, b) \quad \text{for all } 0 \leq a \leq b \leq 1,$$

where $C$ does not depend on $a$ or $b$. Hence

$$E\{\widetilde{K}'_L((\alpha_1, \alpha_2] \times (r_1, r])^2\}$$
$$\leq \frac{C}{L} \sum_{l=[Lr_1]+1}^{[Lr]} \{|\alpha_2 - \alpha_1| + q(\alpha_1, \alpha_2)\}$$
$$= C\frac{[Lr] - [Lr_1]}{L} \{|\alpha_2 - \alpha_1| + q(\alpha_1, \alpha_2)\}.$$

Note that, since we can restrict to $(r_2 - r_1) \geq \frac{1}{L}$,
$$\frac{[Lr] - [Lr_1]}{L} \leq 2(r_2 - r_1).$$
This bound obviously holds for $B_2$ as well, hence
$$E(\widetilde{K}'_L((\alpha_1, \alpha_2] \times (r_1, r])^2 \widetilde{K}'_L((\alpha_1, \alpha_2] \times (r, r_2])^2)$$
$$\leq C(r_2 - r_1)\{|\alpha_2 - \alpha_1| + q(\alpha_1, \alpha_2)\}^2.$$
Let us now consider Type I blocks: we have
$$E\{\widetilde{K}'_L((\alpha_1, \alpha] \times (r_1, r_2])^2 \widetilde{K}'_L((\alpha, \alpha_2] \times (r_1, r_2])^2\}$$
$$= \frac{1}{L^2} \sum_{l=[Lr_1]+1}^{[Lr_2]} EG_l^{*2}((\alpha_1, \alpha] \cap (0, 1 - l^{-3/2}])$$
$$\times G_l^{*2}((\alpha, \alpha_2] \cap (0, 1 - l^{-3/2}])$$
$$+ \left\{ \frac{1}{L} \sum_{l=[Lr_1]+1}^{[Lr_2]} EG_l^{*2}((\alpha_1, \alpha] \cap (0, 1 - l^{-3/2}]) \right\}$$
$$\times \left\{ \frac{1}{L} \sum_{l=[Lr_1]+1}^{[Lr_2]} EG_l^{*2}((\alpha, \alpha_2] \cap (0, 1 - l^{-3/2}]) \right\}$$
$$+ \left\{ \frac{1}{L} \sum_{l=[Lr_1]+1}^{[Lr_2]} EG_l^*((\alpha_1, \alpha] \cap (0, 1 - l^{-3/2}]) \right.$$
$$\left. \times G_l^*((\alpha, \alpha_2] \cap (0, 1 - l^{-3/2}]) \right\}^2.$$

The middle term is the same as for Type II blocks, while the last term is bounded by the former because of the Cauchy–Schwarz inequality. As far as the first term is concerned, notice that
$$EG_l^{*2}((\alpha_1, \alpha] \cap (0, 1 - l^{-3/2}])G_l^{*2}((\alpha, \alpha_2] \cap (0, 1 - l^{-3/2}]) = 0,$$
if $\alpha > 1 - l^{-3/2}$. Otherwise, the right-hand side is equal to
$$EG_l^{*2}(\alpha_1, \alpha]G_l^{*2}(\alpha, \tau_l(\alpha_2)]$$

(18) $$= \frac{1}{l} EZ_{l1}^2(\alpha_1, \alpha] Z_{l1}^2(\alpha, \tau_l(\alpha_2)]$$

(19) $$+ 2\frac{l-1}{l} EZ_{l1}^2(\alpha_1, \alpha] Z_{l1}(\alpha, \tau_l(\alpha_2)) Z_{l2}(\alpha, \tau_l(\alpha_2))$$

(20) $$+ 2\frac{l-1}{l} EZ_{l1}(\alpha_1, \alpha] Z_{l1}^2(\alpha, \tau_l(\alpha_2)) Z_{l2}(\alpha_1, \alpha]$$

$$(21) \quad + \frac{l-1}{l}(l-2)EZ_{l1}^2(\alpha_1,\alpha]Z_{l2}(\alpha,\tau_l(\alpha_2)]Z_{l3}(\alpha,\tau_l(\alpha_2)]$$

$$(22) \quad + \frac{l-1}{l}(l-2)EZ_{l1}^2(\alpha,\tau_l(\alpha_2)]Z_{l2}(\alpha_1,\alpha]Z_{l3}(\alpha_1,\alpha]$$

$$(23) \quad + 4\frac{l-1}{l}EZ_{l1}(\alpha_1,\alpha]Z_{l1}(\alpha,\tau_l(\alpha_2)]Z_{l2}(\alpha_1,\alpha]Z_{l2}(\alpha,\tau_l(\alpha_2)]$$

$$(24) \quad + \frac{l-1}{l}EZ_{l1}^2(\alpha_1,\alpha]Z_{l2}^2(\alpha,\tau_l(\alpha_2)]$$

$$(25) \quad + 4\frac{l-1}{l}(l-2)$$
$$\times EZ_{l1}(\alpha_1,\alpha]Z_{l1}(\alpha,\tau_l(\alpha_2)]Z_{l2}(\alpha_1,\alpha]Z_{l2}(\alpha,\tau_l(\alpha_2)]$$

$$(26) \quad + \frac{l-1}{l}(l-2)(l-3)$$
$$\times EZ_{l1}(\alpha_1,\alpha]Z_{l2}(\alpha_1,\alpha]Z_{l3}(\alpha,\tau_l(\alpha_2)]Z_{l4}(\alpha,\tau_l(\alpha_2)].$$

Now for (18) note that

$$(27) \quad \mathbb{1}(\alpha_1 < y_{l1} \leq \alpha)\mathbb{1}(\alpha < y_{l1} \leq \alpha_2) \equiv 0,$$
$$p_l^2(\alpha,\beta) \leq p_l(\alpha,\beta) \leq 1,$$

and then

$$EZ_{l1}^2(\alpha_1,\alpha]Z_{l1}^2(\alpha,\tau_l(\alpha_2)]$$
$$= E\{\mathbb{1}(\alpha_1 < y_{l1} \leq \alpha) - p_l(\alpha_1,\alpha)\}^2$$
$$\times \{\mathbb{1}(\alpha < y_{l2} \leq \tau_l(\alpha_2)) - p_l(\alpha,\tau_l(\alpha_2))\}^2$$
$$\leq Cp_l(\alpha_1,\alpha)p_l(\alpha,\tau_l(\alpha_2)) \leq C|\alpha-\alpha_1|\,|\alpha_2-\alpha|.$$

For (19), in view of (27), simple manipulations and Lemma A.6,

$$E\{(\mathbb{1}(\alpha_1 < y_{l1} \leq \alpha) - p_l(\alpha_1,\alpha))^2 Z_{l1}(\alpha,\tau_l(\alpha_2)]Z_{l2}(\alpha,\tau_l(\alpha_2)]\}$$
$$= (1 - 2p_l(\alpha_1,\alpha))E\{\mathbb{1}(\alpha_1 < y_{l1} \leq \alpha)Z_{l1}(\alpha,\tau_l(\alpha_2)]Z_{l2}(\alpha,\tau_l(\alpha_2)]\}$$
$$+ p_l(\alpha_1,\alpha)^2 E\{Z_{l1}(\alpha,\tau_l(\alpha_2)]Z_{l2}(\alpha,\tau_l(\alpha_2)]\}$$
$$= -(1 - 2p_l(\alpha_1,\alpha))p_l(\alpha,\alpha_2)E\{\mathbb{1}(\alpha_1 < y_{l1} \leq \alpha)Z_{l2}(\alpha,\tau_l(\alpha_2)]\}$$
$$+ p_l(\alpha_1,\alpha)^2 E\{Z_{l1}(\alpha,\tau_l(\alpha_2)]Z_{l2}(\alpha,\tau_l(\alpha_2)]\}$$
$$\leq Cp_l(\alpha_1,\alpha)q(\alpha,\tau_l(\alpha_2)) \leq C|\alpha-\alpha_1|q(\alpha,\alpha_2).$$

The argument for (20), (23) and (24) is entirely analogous. For (21) we obtain

$$EZ_{l1}^2(\alpha_1,\alpha]Z_{l2}(\alpha,\tau_l(\alpha_2)]Z_{l3}(\alpha,\tau_l(\alpha_2)]$$

$$= (1 - p_l(\alpha_1, \alpha))E\{\mathbb{1}(\alpha_1 < y_{l1} \le \alpha)Z_{l1}(\alpha, \tau_l(\alpha_2)]Z_{l2}(\alpha, \tau_l(\alpha_2)]\}$$
$$- p_l(\alpha_1, \alpha)E\{Z_{l1}(\alpha_1, \alpha]Z_{l1}(\alpha, \tau_l(\alpha_2)]Z_{l2}(\alpha, \tau_l(\alpha_2)]\}$$
$$= (1 - 2p_l(\alpha_1, \alpha))E\{Z_{l1}(\alpha_1, \alpha]Z_{l1}(\alpha, \tau_l(\alpha_2)]Z_{l2}(\alpha, \tau_l(\alpha_2)]\}$$
$$+ (1 - p_l(\alpha_1, \alpha))p_l(\alpha_1, \alpha)E\{Z_{l1}(\alpha, \tau_l(\alpha_2)]Z_{l2}(\alpha, \tau_l(\alpha_2)]\}$$
$$\le \frac{C}{l}p_l(\alpha_1, \alpha)q(\alpha_1, \alpha)q(\alpha, \alpha_2) \le \frac{C}{l}|\alpha - \alpha_1|q(\alpha_1, \alpha)q(\alpha, \alpha_2),$$

in view of previous results and Lemma A.7. The argument for (22) and (25) is analogous. Finally, it is clear that for (26) it is sufficient to prove

$$EZ_{l1}(\alpha_1, \alpha]Z_{l2}(\alpha_1, \alpha]Z_{l3}(\alpha, \tau_l(\alpha_2)]Z_{l4}(\alpha, \tau_l(\alpha_2)]$$
$$\le \frac{C}{l^2}q^2(\alpha_1, \alpha)q^2(\alpha, \tau_l(\alpha_2)) \le \frac{C}{l^2}q(\alpha_1, \alpha)q(\alpha, \alpha_2),$$

a result established in Lemma A.7. Hence each of the terms in (18)–(26) is bounded uniformly by

$$C\{|\alpha_2 - \alpha_1| + q(\alpha_1, \alpha_2)\}^2.$$

The proof can then be completed by routine manipulations. □

REMARK 3.1. It is immediately seen that the bounds of (18)–(26) given previously hold for arbitrary nonoverlapping blocks. As such they can be used to establish the Lindeberg condition (13) in the proof of Proposition 3.1.

**4. Comments and conclusions.** Theorem 2.1 is immediately applicable to statistical inference procedures, and in particular, to tests for Gaussianity. For instance, a Kolmogorov–Smirnov type of test is implemented, for any suitably large $L$, if we evaluate

(28) $$S_L = \sup_{\alpha, r} |K_L^*(\alpha, r)|,$$

and compare the observed value with the desired quantile of the law of $S_\infty = \sup_{\alpha, r} |K^*(\alpha, r)|$. In principle, the latter value can be derived by Monte Carlo simulation, since the limiting distribution does not entail any unknown nuisance parameter. Likewise, Cramer–Von Mises and many other types of goodness-of-fit statistics could be easily implemented.

Tests for Gaussianity on spherical maps have recently been considered by several authors in the physics literature, as mentioned in the Introduction. The focus of these papers is very much on the physical discussion rather than on statistical methodology, so that any comparison seems inappropriate. We claim, however, that our method enjoys some important advantages on any empirical procedure in this literature. To mention a few, our procedure allows for a rigorous asymptotic theory, which is made possible by the focus on harmonic coefficients; it allows for inference completely free from nuisance parameters, whereas in other papers test statistics are considered whose law depends on the values of the angular power spectrum $C_l$. Due to our study of the asymptotic behavior of the whole field $K_L^*$, many different testing procedures can be implemented; these procedures provide information not only on departures from Gaussianity, but also on their location in harmonic space; this is important, as different physical mechanisms are known to operate at the various multipoles. Moreover, the distributional properties of the normalized spherical harmonic coefficients may have some independent interest for other cosmological applications.

It is clearly of great interest to evaluate the power of our testing procedures in non-Gaussian situations. With regard to this issue, a crucial point is the nature of non-Gaussianity. In the cosmological literature departures from Gaussianity are occasionally generated by superimposing non-Gaussian structures over a Gaussian map. For instance, it is possible to mimic a popular class of topological defects models (the so-called "cosmic strings") by setting

$$(29) \qquad T(\theta, \varphi) = T^{\mathrm{G}}(\theta, \varphi) + T^{\mathrm{S}}(\theta, \varphi),$$

where $T^{\mathrm{G}}(\theta, \varphi)$ is a Gaussian map, and $T^{\mathrm{S}}(\theta, \varphi)$ is a map of Poisson-distributed segments of randomly varying directions, length and level [for more details on simulations of non-Gaussian spherical fields and their physical meaning, see Hansen, Marinucci, Natoli and Vittorio (2002)]. $T^{\mathrm{G}}(\theta, \varphi)$ and $T^{\mathrm{S}}(\theta, \varphi)$ are taken to be zero-mean and independent; we define the percentage of non-Gaussianity in the map by

$$P_{\mathrm{NG}} = \frac{ET^{\mathrm{S}}(\theta, \varphi)^2}{ET^{\mathrm{G}}(\theta, \varphi)^2 + ET^{\mathrm{S}}(\theta, \varphi)^2}.$$

We fix $L = 500$, a realistic value for the MAP experiment, and we evaluate the threshold values for sizes $\alpha = 10\%$, $5\%$, $1\%$ by 500 Monte Carlo replications of (28): we obtain 0.947, 1.012, 1.160, respectively. We then generate string maps with 100, 500 and 1000 strings, and report the rejection frequencies for 100 replications of model (29) (with different percentages of non-Gaussianity) in Table 1.

We leave several issues for future research. To increase the power of the tests, it is useful to consider empirical processes constructed on several rows [see the applied papers by Hansen, Marinucci, Natoli and Vittorio (2002) and Hansen, Marinucci and Vittorio (2003)]. As these rows are independent under the null of Gaussianity, the extension is conceptually straightforward, although computationally burdensome. For cosmological applications, some additional difficulties may arise in practical implementation; in particular, it may be the case that $T(\theta, \varphi)$ is subject to some measurement error, or that it is only incompletely observed, or both. As mentioned before, with the current status of technology on satellite- and balloon-borne experiments, it is known that observational error is some order of magnitude smaller than signal, and thus can safely be ignored, for $l$ as large as 1000 or more; we hence consider the asymptotic theory provided in this paper to be valuable for the practitioner, as $L$ on the order of $10^3$ (with a total number of observed $a_{lm}$ on the order of $10^5/10^6$) seems sufficient for the asymptotic theory to yield applicable approximations. However, the approximation of these sampling distributions can possibly be improved by bootstrap methods, to take into account the nuisance parameters governing noise. The presence of gaps in the observed field poses more challenging questions, and will be addressed elsewhere.

TABLE 1
*Power of the Kolmogorov–Smirnov test*

| Size of the Test\$P_{\text{NG}}$ | 0.1 | 0.2 | 0.3 | 0.4 | 0.5 |
|---|---|---|---|---|---|
| 100 cosmic strings | | | | | |
| 1% | 28% | 83% | 98% | 100% | 100% |
| 5% | 38% | 86% | 99% | 100% | 100% |
| 10% | 42% | 88% | 100% | 100% | 100% |
| 500 cosmic strings | | | | | |
| 1% | 6% | 41% | 82% | 92% | 100% |
| 5% | 14% | 49% | 88% | 97% | 100% |
| 10% | 20% | 56% | 90% | 99% | 100% |
| 1000 cosmic strings | | | | | |
| 1% | 3% | 23% | 66% | 95% | 100% |
| 5% | 8% | 35% | 78% | 97% | 100% |
| 10% | 16% | 40% | 81% | 97% | 100% |

## APPENDIX

LEMMA A.1. *Let $a_n$ be a sequence of real numbers such that $\lim_{n\to\infty} n^2 a_n = -c < \infty$. Then*

$$\lim_{n\to\infty} n((1+a_n)^n - 1) = c.$$

PROOF. The proof is trivial, writing $(1+a_n)^n = \exp(n\log(1+a_n))$ and using the expansions, as $r \to 0$

(30) $$\log(1+r) = r - \frac{r^2}{2} + O(r^3),$$

(31) $$e^r = 1 + r + O(r^2).$$

□

The proof of the following lemma was very much shortened thanks to the comments of one referee.

LEMMA A.2. *For any $C > 0$, as $\to \infty$,*

$$\sup_{l^{-C} \leq x \leq 1} \left| l\left\{ x - \left(1 + \frac{\log x}{l}\right) - x\log x - \frac{x\log^2 x}{2}\right\} \right| \to 0.$$

PROOF. We have, uniformly in $l^{-C} \leq x \leq 1$,

$$l\left\{x - \left(1 + \frac{\log x}{l}\right)\right\} = l\left\{x - \exp\left((l-1)\log\left(1 + \frac{\log x}{l}\right)\right)\right\}$$

$$= l\left\{x - \exp\left((l-1)\left(\frac{\log x}{l} - \frac{\log^2 x}{2l^2} + O\left(\frac{\log^3 l}{l^3}\right)\right)\right)\right\}$$

$$= l\left\{x - \exp\left(\log x - \frac{\log x}{l} - \frac{\log^2 x}{2l} + O\left(\frac{\log^3 l}{l^2}\right)\right)\right\}$$

$$= lx\left\{1 - \exp\left(-\frac{\log x}{l} - \frac{\log^2 x}{2l} + O\left(\frac{\log^3 l}{l^2}\right)\right)\right\}$$

$$= lx\left\{1 - \left(1 - \frac{\log x}{l} - \frac{\log^2 x}{2l} + O\left(\frac{\log^4 l}{l^2}\right)\right)\right\}$$

$$= x\log x + \frac{x\log^2 x}{2} + O\left(\frac{\log^4 l}{l}\right)$$

$$= x\log x + \frac{x\log^2 x}{2} + o(1),$$

where we used (30) and (31). □

LEMMA A.3. *The marginal density of $y_{lm} = 1 - \exp(-l\xi_{lm})$ is bounded uniformly by $e^2$ for $l \geq 2$.*

PROOF. The marginal density of $\xi_{lm}$ is

$$f_{\xi_{lm}}(t) = (l-1)(1-t)^{l-2}\mathbb{1}_{[0,1]}(t),$$

whence by the change of variable formula,

$$f_{y_{lm}}(y) = \frac{l-1}{l(1-y)}\left(1 + \frac{\log(1-y)}{l}\right)^{l-2}\mathbb{1}_{[0,1-\exp(-l))}(y)$$

$$\leq \frac{1}{(1-y)}\left(1 + \frac{\log(1-y)}{l}\right)^{l-2}\mathbb{1}_{[0,1-\exp(-l))}(y).$$

Now let $x = 1 - y$ and differentiate with respect to $x$; we obtain

$$\frac{(1 - 2/l)(1 + l^{-1} \log x)^{l-3} - (1 + l^{-1} \log x)^{l-2}}{x^2}$$
$$= -\frac{(1 + l^{-1} \log x)^{l-3} l^{-1} (2 + \log x)}{x^2},$$

which is equal to zero at $x = e^{-2}$; the latter is easily seen to be a unique maximum. Hence

$$f_{Y_l}(y) \leq e^2 \left(1 + \frac{\log e^{-2}}{l}\right)^{l-2} \leq e^2. \qquad \square$$

LEMMA A.4. *As $L \to \infty$,*

$$\sup_{\alpha, r} |\widehat{K}_L(\alpha, r) - \widetilde{K}_L(\alpha, r)| = o_p(1).$$

PROOF. By the definition of $\widetilde{K}_L(\alpha, r)$ (14),

$$\sup_{\alpha, r} |\widehat{K}_L(\alpha, r) - \widetilde{K}_L(\alpha, r)| \leq \frac{1}{\sqrt{L}} \sum_{l=1}^{L} \frac{1}{\sqrt{l}} \sum_{m=1}^{l} \mathbb{1}(1 - l^{-3/2} < y_{lm} \leq 1);$$

hence, using Lemma A.3 and Markov's inequality,

$$\sup_{\alpha, r} |\widehat{K}_L(\alpha, r) - \widetilde{K}_L(\alpha, r)| = O_p\left(\sum_{l=1}^{L} \sqrt{l} E\{\mathbb{1}(0 \leq e^{-l\xi_{lm}} < l^{-3/2})\}\right)$$
$$= O_p\left(\frac{1}{\sqrt{L}} \sum_{l=1}^{L} \frac{1}{l}\right) = O_p\left(\frac{\log L}{\sqrt{L}}\right), \quad \text{as } L \to \infty.$$

$\square$

In the sequel, recall that

$$Z_{lm}((\alpha, \beta]) = \mathbb{1}(\alpha_1 < y_{lm} \leq \beta) - p_l(\alpha, \beta),$$
$$p_l(\alpha, \beta) = E\mathbb{1}(\alpha_1 < y_{lm} \leq \beta).$$

LEMMA A.5. *As $L \to \infty$,*

$$\sup_{\alpha, r} |E\widetilde{K}_L(\alpha, r) - 2\sqrt{r} b(\alpha)| = o(1).$$

PROOF. By the definition of $\widetilde{K}_L(\alpha, r)$,

$$E\widetilde{K}_L(\alpha, r) - 2\sqrt{r}b(\alpha)$$

$$= \frac{1}{\sqrt{L}} \sum_{l=1}^{[Lr]} \left[\frac{1}{\sqrt{l}} b_l(t_l(\alpha))\right] - 2\sqrt{r}b(\alpha)$$

$$= \frac{1}{\sqrt{L}} \sum_{l=1}^{[Lr]} \frac{1}{\sqrt{l}} [b_l(t_l(\alpha)) - b(\alpha)] + b(\alpha) \left[\frac{1}{\sqrt{L}} \sum_{l=1}^{[Lr]} \frac{1}{\sqrt{l}} - 2\sqrt{r}\right],$$

whose absolute value is bounded by

$$(32) \qquad \frac{1}{\sqrt{L}} \sum_{l=1}^{L} \frac{1}{\sqrt{l}} \sup_{\alpha} |b_l(t_l(\alpha)) - b(\alpha)|$$

$$(33) \qquad + \sup_{\alpha} |b(\alpha)| \sup_{r} \left|\frac{1}{\sqrt{L}} \sum_{l=1}^{[Lr]} \frac{1}{\sqrt{l}} - 2\sqrt{r}\right|.$$

Now for (33) we have that $\sup_\alpha |b(\alpha)| \le C$, whereas approximating the sum with the integral, it is easy to see that

$$\frac{1}{\sqrt{L}} \sum_{l=1}^{[Lr]} \frac{1}{\sqrt{l}} = \frac{1}{L} \sum_{l=1}^{[Lr]} \frac{1}{\sqrt{l/L}} \le \int_0^r \frac{1}{\sqrt{x}} dx = 2\sqrt{r},$$

and

$$2\sqrt{r} = \int_0^{1/L} \frac{1}{\sqrt{x}} dx + \int_{1/L}^r \frac{1}{\sqrt{x}} dx$$

$$\le \frac{2}{\sqrt{L}} + \frac{1}{L} \sum_{l=2}^{[Lr]} \frac{1}{\sqrt{(l-1)/L}}.$$

Hence

$$(33) \le C \frac{2}{\sqrt{L}} = o(1) \qquad \text{as } L \to \infty.$$

On the other hand, for (32) it is enough to prove that

$$\sup_{0 \le \alpha \le 1} |b_l(t_l(\alpha)) - b(\alpha)|$$

$$\le \sup_{0 \le \alpha \le 1-l^{-3/2}} |b_l(t_l(\alpha)) - b(\alpha)| + \sup_{1-l^{-3/2} \le \alpha \le 1} |b(\alpha)|$$

is $o(1)$ as $l \to \infty$. For the second term on the right-hand side, continuity of $b(\alpha)$ at $\alpha = 1$ is enough. For the first, we take $x = 1 - \alpha$, so that we need to establish that

$$\lim_{l \to \infty} \sup_{l^{-3/2} \le x \le 1} \left| l\left\{ x - \left(1 + \frac{\log x}{l}\right)^{l-1}\right\} - \frac{1}{2} x \log^2 x - x \log x \right| = 0,$$

which follows from Lemma A.2. Thus the proof is completed. □

LEMMA A.6. *For $0 \leq \alpha_1 < \alpha_2 \leq 1 - l^{-3/2}$, $0 \leq \beta_1 < \beta_2 \leq 1 - l^{-3/2}$, we have*

$$|\mathrm{Cov}\{\mathbb{1}(\alpha_1 \leq y_{l1} \leq \alpha_2), \mathbb{1}(\beta_1 \leq y_{l2} \leq \beta_2)\}| \leq \frac{C}{l} q(\alpha_1, \alpha_2) q(\beta_1, \beta_2),$$

*where*

$$q(a,b) = \int_a^b \{1 + |\log(1-y)|\}\, dy$$

$$= \int_{1-b}^{1-a} \{1 + |\log x|\}\, dx$$

*is a finite measure and $C$ does not depend on $l, \alpha, \beta$.*

PROOF. After a change of variable we have

$$\mathrm{Cov}\{\mathbb{1}(\alpha_1 \leq y_{l1} \leq \alpha_2), \mathbb{1}(\beta_1 \leq y_{l2} \leq \beta_2)\} = \int_{1-\alpha_2}^{1-\alpha_1} \int_{1-\beta_2}^{1-\beta_1} h_l(x_1, x_2)\, dx_1\, dx_2,$$

where the covariance density $h_l(x_1, x_2)$ is given by

$$\frac{1}{x_1 x_2} \left\{ \frac{(l-1)(l-2)}{l^2} \left(1 + \frac{\log x_1 x_2}{l}\right)^{l-3} \right.$$
$$\left. - \frac{(l-1)^2}{l^2} \left(1 + \frac{\log x_1}{l}\right)^{l-2} \left(1 + \frac{\log x_2}{l}\right)^{l-2} \right\}$$
$$= A_l(x_1, x_2) \times B_l(x_1, x_2),$$

for

$$A_l(x_1, x_2) = \left\{ \left(1 + \frac{\log x_1 x_2}{l}\right)^3 \left(1 + \frac{\log x_1}{l}\right)^2 \left(1 + \frac{\log x_2}{l}\right)^2 \right\}^{-1}$$
$$\times \frac{1}{x_1 x_2} \left(1 + \frac{\log x_1 x_2}{l}\right)^l,$$

$$B_l(x_1, x_2) = \left\{ \frac{(l-1)(l-2)}{l^2} \left(1 + \frac{\log x_1}{l}\right)^2 \left(1 + \frac{\log x_2}{l}\right)^2 \right.$$
$$\left. - \frac{(l-1)^2}{l^2} \frac{(1 + l^{-1}\log x_1)^2 (1 + l^{-1}\log x_2)^2}{(1 + l^{-1}\log x_1 x_2)^l} \right\},$$

$$l^{-3/2} \leq x_1, x_2 \leq 1.$$

Now for $A_l(x_1, x_2)$ we have

$$\left(1 + \frac{\log x_1 x_2}{l}\right)^l \leq x_1 x_2,$$

because

(34) $$1 + \frac{\log x}{l} \leq \exp\left(\frac{\log x}{l}\right),$$

always. Also

$$\sup_{l^{-3/2} \leq x_1, x_2 \leq 1} \left\{\left(1 + \frac{\log x_1 x_2}{l}\right)^3 \left(1 + \frac{\log x_1}{l}\right)^2 \left(1 + \frac{\log x_2}{l}\right)^2\right\}^{-1}$$

$$\leq \left\{\left(1 - \frac{9}{4}\frac{\log l}{l}\right)^3 \left(1 - \frac{3}{2}\frac{\log l}{l}\right)^4\right\}^{-1} \leq C,$$

uniformly for $l \geq 1$. Thus $A_l(x_1, x_2)$ is bounded uniformly in $x$ and $l$. Now call

$$c_{1,l} = \frac{(l-1)(l-2)}{l^2}, \qquad c_{2,l} = \frac{(l-1)^2}{l}.$$

Thus $B_l(x_1, x_2)$ becomes

(35)
$$c_{2,l}\left(1 + 2\frac{\log x_1}{l} + \frac{\log^2 x_1}{l^2}\right)\left(1 + 2\frac{\log x_2}{l} + \frac{\log^2 x_2}{l^2}\right)$$
$$- c_{1,l}\left(1 + \frac{\log x_1 \log x_2}{l^2 + l \log x_1 x_2}\right)^l$$

and because $|c_{2,l} - c_{1,l}| \leq l^{-1}$,

(36) $\sup_{l^{-3/2} \leq x \leq 1} \frac{|\log x|}{l} \leq 1, \qquad \frac{\log^2 x}{l^2} \leq \frac{|\log x|}{l} \qquad$ for all $l^{-3/2} \leq x \leq 1$,

$$\left(1 + \frac{\log x_1 \log x_2}{l^2 + l \log x_1 x_2}\right)^l = 1 + \sum_{k=1}^{l} \binom{l}{k} \left(\frac{\log x_1 \log x_2}{l^2 + l \log x_1 x_2}\right)^k$$

$$= 1 + \frac{\log x_1 \log x_2}{l + \log x_1 x_2} + \sum_{k=2}^{l} \binom{l}{k} \left(\frac{\log x_1 \log x_2}{l^2 + l \log x_1 x_2}\right)^k,$$

we obtain

$$|(35)| \leq |c_{2,l} - c_{1,l}| + \frac{C}{l}\{|\log x_1| + |\log x_2| + |\log x_1||\log x_2|\}$$

$$+ c_{1,l}\left|\sum_{k=1}^{l} \binom{l}{k}\left(\frac{\log x_1 \log x_2}{l^2 + l \log x_1 x_2}\right)^k\right|$$

$$\leq \frac{C}{l}\{1 + |\log x_1| + |\log x_2| + |\log x_1||\log x_2|\}$$

$$+ C\frac{|\log x_1||\log x_2|}{l - 3\log l}\left\{1 + \frac{1}{l}\left|\sum_{k=2}^{l} \binom{l}{k}\left(\frac{|\log x_1 \log x_2|}{l^2 + l \log x_1 x_2}\right)^{k-1}\right|\right\}$$

$$\leq \frac{C}{l}\{1 + |\log x_1| + |\log x_2| + |\log x_1||\log x_2|\}.$$

The last step follows from

$$\sup_{l^{-3/2} \leq x_1, x_2 \leq 1} \frac{1}{l} \sum_{k=2}^{l} \binom{l}{k} \left( \frac{|\log x_1 \log x_2|}{l^2 + l \log x_1 x_2} \right)^{k-1}$$

$$\leq \sup_{l^{-3/2} \leq x_1, x_2 \leq 1} \sum_{k=2}^{l} \frac{l^{k-1}}{k!} \left( \frac{|\log x_1 \log x_2|}{l^2 + l \log x_1 x_2} \right)^{k-1}$$

$$= \sum_{k=2}^{l} \frac{1}{k!} \left( \frac{(9/4) \log^2 l}{l - 3 \log l} \right)^{k-1} \leq C,$$

because $\log^2 l/(l - 3 \log l)$ is uniformly bounded in $l$.  $\square$

LEMMA A.7. *For all $0 \leq \alpha_1 < \alpha_2 \leq 1 - l^{-3/2}$, $0 \leq \beta_1 < \beta_2 \leq 1 - l^{-3/2}$, $0 \leq \gamma_1 < \gamma_2 \leq 1 - l^{-3/2}$, $0 \leq \delta_1 < \delta_2 \leq 1 - l^{-3/2}$, we have*

(37)
$$|EZ_{l1}(\alpha_1, \alpha_2]Z_{l2}(\beta_1, \beta]Z_{l3}(\gamma_1, \gamma_2]|$$
$$\leq \frac{C}{l} q(\alpha_1, \alpha_2) q(\beta_1, \beta) q(\gamma_1, \gamma_2),$$

*and*

(38)
$$|EZ_{l1}(\alpha_1, \alpha_2]Z_{l2}(\beta_1, \beta]Z_{l3}(\gamma_1, \gamma_2]Z_{l4}(\delta_1, \delta_2]|$$
$$\leq \frac{C}{l^2} q(\alpha_1, \alpha_2) q(\beta_1, \beta) q(\gamma_1, \gamma_2) q(\delta_1, \delta_2),$$

*where $C$ does not depend on $\alpha$, $\beta$, $\gamma$, $\delta$ or $l$.*

PROOF. For the sake of brevity, we give only the proof of (38); the proof of (37) is similar, indeed slightly simpler, and can be found in Marinucci and Piccioni (2003). After a change of variable we have

$$EZ_{l1}(\alpha_1, \alpha_2]Z_{l2}(\beta_1, \beta]Z_{l3}(\gamma_1, \gamma_2]Z_{l4}(\delta_1, \delta_2]$$
$$= \int_{1-\alpha_2}^{1-\alpha_1} \int_{1-\beta_2}^{1-\beta_1} \int_{1-\gamma_2}^{1-\gamma_1} \int_{1-\delta_2}^{1-\delta_1} h_l(x_1, x_2, x_3, x_4) \, dx_1 \, dx_2 \, dx_3 \, dx_4,$$

where the four-term covariance density $h_l(x_1, x_2, x_3, x_4)$ is given by, for $l^{-3/2} \leq x_1, x_2, x_3, x_4 \leq 1$,

$$\frac{1}{x_1 x_2 x_3 x_4} \frac{(l-1)(l-2)(l-3)(l-4)}{l^4} \left(1 + \frac{\log x_1 x_2 x_3 x_4}{l}\right)^{l-5}$$

$$- \frac{1}{x_1 x_2 x_3 x_4} \sum_{i=1}^{4} \frac{(l-1)^2(l-2)(l-3)}{l^4}$$

$$\times \left(1 + \frac{\log x_i}{l}\right)^{l-2} \left(1 + \frac{\sum_{j\neq i} \log x_j}{l}\right)^{l-4}$$

$$+ \frac{1}{x_1 x_2 x_3 x_4} \sum_{i=1}^{3} \sum_{j=i+1}^{4} \frac{(l-1)^3(l-2)}{l^4} \left(1 + \frac{\log x_i}{l}\right)^{l-2}$$

$$\times \left(1 + \frac{\log x_j}{l}\right)^{l-2} \left(1 + \frac{\sum_{k\neq j,i} \log x_k}{l}\right)^{l-4}$$

$$- 3\frac{1}{x_1 x_2 x_3 x_4} \frac{(l-1)^4}{l^4} \left(1 + \frac{\log x_1}{l}\right)^{l-2} \left(1 + \frac{\log x_2}{l}\right)^{l-2}$$

$$\times \left(1 + \frac{\log x_3}{l}\right)^{l-2} \left(1 + \frac{\log x_4}{l}\right)^{l-2}.$$

Define

$$c_{1l} = \frac{(l-1)^4}{l^4}, \qquad c_{2l} = \frac{(l-1)^3(l-2)}{l^4},$$

$$c_{3l} = \frac{(l-1)^2(l-2)(l-3)}{l^4}, \qquad c_{4l} = \frac{(l-1)(l-2)(l-3)(l-4)}{l^4}.$$

We have

$$h_l(x_1, x_2, x_3, x_4)$$
$$= \psi_l \left\{ c_{4l} \zeta_l - \sum_{i=1}^{4} c_{3l} \eta_{il} d_{il} + \sum_{i=1}^{3} \sum_{j=i+1}^{4} c_{2l} \theta_{ijl} d_{ijl} - 3 c_{4l} \xi_l d_l \right\},$$

where

$$\psi_l = \zeta_l^{-1} a_l (x_1 x_2 x_3 x_4)^{l-5}, \qquad a_l(x) = 1 + \frac{\log x}{l},$$

$$\zeta_l = \left(\prod_{i=1}^{4} a_l^2(x_i)\right) \left(\prod_{i=1}^{3} \prod_{j=i+1}^{4} a_l^3(x_i x_j)\right) \left(\prod_{i=1}^{2} \prod_{j=i+1}^{3} \prod_{k=j+1}^{4} a_l^4(x_i x_j x_k)\right),$$

$$\eta_{il} = \left(\prod_{j\neq i}^{4} a_l^2(x_j)\right) \left(\prod_{j=1}^{3} \prod_{k=j+1}^{4} a_l^3(x_i x_j)\right) \left(\prod_{j,k\neq i} a_l^4(x_i x_j x_k)\right) a_l^5(x_1 x_2 x_3 x_4),$$

$$\theta_{ijl} = \left(\prod_{k\neq i,j} a_l^2(x_k)\right) \left(\prod_{k,l\neq i,j} a_l^3(x_k x_l)\right)$$
$$\times \left(\prod_{k=1}^{2} \prod_{l=k+1}^{3} \prod_{m=l+1}^{4} a_l^4(x_k x_l x_m)\right) a_l^5(x_1 x_2 x_3 x_4),$$

$$\xi_l = \left( \prod_{i=1}^{3} \prod_{j=i+1}^{4} a_l^3(x_i x_j) \right) \left( \prod_{i=1}^{2} \prod_{j=i+1}^{3} \prod_{k=j+1}^{4} a_l^4(x_i x_j x_k) \right) a_l^5(x_1 x_2 x_3 x_4),$$

and

$$d_{il} = \left[ \frac{(l + \log x_i)(l + \log \prod_{j \neq i} x_j)}{l^2 + l \log x_1 x_2 x_3 x_4} \right]^l,$$

$$d_{ijl} = \left[ \frac{(l + \log x_i)(l + \log x_j)(l + \log \prod_{k \neq i,j} x_k)}{l^3 + l^2 \log x_1 x_2 x_3 x_4} \right]^l,$$

$$d_l = \left[ \frac{(l + \log x_1)(l + \log x_2)(l + \log x_3)(l + \log x_4)}{l^4 + l^3 \log x_1 x_2 x_3 x_4} \right]^l.$$

Now note that

$$d_{il} = 1 + \sum_{k=1}^{l} \binom{l}{k} \left[ \frac{\log x_i (\sum_{j \neq i} \log x_j)}{l^2 + l \log x_1 x_2 x_3 x_4} \right]^k$$

$$= 1 + \frac{\log x_i (\sum_{j \neq i} \log x_j)}{l + \log x_1 x_2 x_3 x_4} + e_l\left( \log x_i; \sum_{j \neq i} \log x_j \right),$$

where

$$\left| e_l\left( \log x_i; \sum_{j \neq i} \log x_j \right) \right|$$

$$= \left| \sum_{k=2}^{l} \binom{l}{k} \left[ \frac{\log x_i (\sum_{j \neq i} \log x_j)}{l^2 + l \log x_1 x_2 x_3 x_4} \right]^k \right|$$

$$\leq \frac{C}{l^2} \left| \frac{l(l-1)}{l^2} \left[ \frac{\log x_i (\sum_{j \neq i} \log x_j)}{1 - l^{-1} 6 \log l} \right]^2 \right| \exp\left\{ \frac{27}{4} \frac{\log^2 l}{l - 6 \log l} \right\}$$

$$\leq \frac{C}{l^2} \left( \sum_{i+1}^{3} \sum_{j=i+1}^{4} |\log x_i| |\log x_j| \right)^2.$$

Also

$$a_l^2(x_i) = 1 + 2\frac{\log x_i}{l} + \frac{\log^2 x_i}{l^2},$$

$$a_l^3(x_i x_j) = 1 + 3\frac{\log x_i}{l} + 3\frac{\log x_j}{l} + O\left( \frac{\log^2 x_i + \log^2 x_j}{l^2} \right),$$

$$a_l^4(x_i x_j x_k) = 1 + 4\frac{\log x_i}{l} + 4\frac{\log x_j}{l} + 4\frac{\log x_k}{l}$$

$$+ O\left( \frac{\log^2 x_i + \log^2 x_j + \log^2 x_k}{l^2} \right),$$

$$a_l^5(x_1 x_2 x_3 x_4) = 1 + \frac{5}{l}\sum_{i=1}^{4} \log x_i + O\Big(\frac{\sum_{i=1}^{4} \log^2 x_i}{l^2}\Big),$$

where the $O(\cdot)$ bounds on the remainders are uniform over $x_i$. Counting terms then gives

$$\zeta_l = 1 + (2 + 3 \times 3 + 3 \times 4)\sum_{i=1}^{4} \frac{\log x_i}{l} + \rho_{\zeta l}(x_1, x_2, x_3, x_4)$$

$$= 1 + 23 \sum_{i=1}^{4} \frac{\log x_i}{l} + \rho_{\zeta l}(x_1, x_2, x_3, x_4),$$

$$|\rho_{\zeta l}(x_1, x_2, x_3, x_4)| \le \frac{C}{l^2} \sum_{i=1}^{3} \sum_{j=i+1}^{4} |\log x_i| |\log x_j|.$$

$$\eta_{il} = 1 + (3 \times 4 + 3 \times 3 + 5)\frac{\log x_i}{l}$$

$$+ (2 + 2 \times 4 + 3 \times 3 + 5)\sum_{j \ne i} \frac{\log x_j}{l} + \rho_{\eta l}(x_i; x_j, x_k, x_l)$$

$$= 1 + 26\frac{\log x_i}{l} + 24 \sum_{j \ne i} \frac{\log x_j}{l} + \rho_{\eta l}(x_i; x_j, x_k, x_l),$$

$$\rho_{\eta l}(x_i; x_j, x_k, x_l) \le \frac{C}{l^2} \sum_{i=1}^{3} \sum_{j=i+1}^{4} |\log x_i| |\log x_j|,$$

whence

$$\sum_{i=1}^{4} \eta_{il} = 4 + 98 \sum_{i=1}^{4} \frac{\log x_i}{l} + O\Big(\frac{1}{l^2} \sum_{i=1}^{3} \sum_{j=i+1}^{4} |\log x_i| |\log x_j|\Big).$$

Likewise, for $k, l \ne i, j$,

$$\theta_{ijl} = 1 + 26 \frac{\log x_i + \log x_j}{l} + 25 \frac{\log x_k + \log x_l}{l}$$

$$+ O\Big(\frac{1}{l^2} \sum_{i=1}^{3} \sum_{j=i+1}^{4} |\log x_i| |\log x_j|\Big),$$

$$\sum_{i=1}^{3} \sum_{j=i+1}^{4} \theta_{ijl} = 6 + 153 \sum_{i=1}^{4} \frac{\log x_i}{l}$$

$$+ O\Big(\frac{1}{l^2} \sum_{i=1}^{3} \sum_{j=i+1}^{4} |\log x_i| |\log x_j|\Big).$$

Finally,

$$\xi_l = 1 + 26 \sum_{i=1}^{4} \frac{\log x_i}{l} + O\left(\frac{1}{l^2} \sum_{i=1}^{3} \sum_{j=i+1}^{4} |\log x_i| |\log x_j|\right).$$

Now combining all terms, we obtain that the covariance density is bounded uniformly in absolute value by

$$|c_{4l} - 4c_{3l} + 6c_{2l} - 3c_{1l}| + |-2c_{3l} + 5c_{2l} - 3c_{1l}| \sum_{i=1}^{4} \frac{|\log x_i|}{l}$$

$$+ |23c_{4l} - 98c_{3l} + 153c_{2l} - 78c_{1l}| \sum_{i=1}^{4} \frac{|\log x_i|}{l}$$

$$+ O\left(\frac{1}{l^2} \sum_{i=1}^{3} \sum_{j=i+1}^{4} |\log x_i| |\log x_j|\right).$$

Now

$$c_{4l} - 4c_{3l} + 6c_{2l} - 3c_{1l}$$
$$= \frac{(l-1)(l-2)(l-3)(l-4) - 4(l-1)^2(l-2)(l-3)}{l^4}$$
$$\quad + \frac{6(l-1)^3(l-2) - 3(l-1)^4}{l^4}$$
$$= \frac{35l^2 - 4(17)l^2 + 6 \times 9 - 3 \times 6}{l^4}$$
$$= \frac{3}{l^2} + O(l^{-3}).$$

Similarly,

$$\lim_{l \to \infty} \{-2c_{3l} + 5c_{2l} - 3c_{1l}\}$$
$$= \lim_{l \to \infty} \frac{-2(l-1)^2(l-2)(l-3) + 5(l-1)^3(l-2) - 3(l-1)^4}{l^3}$$
$$= \frac{1}{l} + O(l^{-2}),$$
$$\lim_{l \to \infty} |23c_{4l} - 98c_{3l} + 153c_{2l} - 78c_{1l}|$$
$$= \frac{3}{l^2} + O(l^{-3}).$$

It can thus be concluded that the four-term covariance density is bounded uniformly by

$$h_l(x_1, x_2, x_3, x_4) \leq \frac{C}{l^2}\left\{1 + \sum_{i=1}^{4} |\log x_i| + \sum_{i=1}^{3}\sum_{j=i+1}^{4} |\log x_i|\,|\log x_j|\right\}$$

$$\leq \frac{C}{l^2}\prod_{i=1}^{4}\left\{1 + \sum_{i=1}^{4}|\log x_i|\right\},$$

and thus the proof is completed. $\square$

**Acknowledgments.** We are very grateful to two anonymous referees whose comments greatly improved the presentation of this paper. We are also grateful to F. Hansen, P. Natoli and N. Vittorio, from the Cosmology Group at the University of Rome "Tor Vergata," for bringing this problem to our attention and for many fruitful discussions. Finally, we thank F. Hansen for carrying out the simulations in Section 4.

## REFERENCES

BICKEL, P. J. and WICHURA, M. J. (1971). Convergence criteria for multiparameter stochastic processes and some applications. *Ann. Math. Statist.* **42** 1656–1670. MR383482

DORÈ, O., COLOMBI, S. and BOUCHET, F. R. (2003). Probing cosmic microwave background non-Gaussianity using local curvature. *Monthly Notices R. Astronom. Soc.* **344** 905–916. Available at arxiv.org as astro-ph/0202135.

DUDLEY, R. M. (1999). *Uniform Central Limit Theorems.* Cambridge Univ. Press. MR1720712

HANNAN, E. J. (1970). *Multiple Time Series.* Wiley, New York. MR279952

HANSEN, F. K., MARINUCCI, D., NATOLI, P. and VITTORIO, N. (2002). Testing for non-Gaussianity of the cosmic microwave background in harmonic space: An empirical process approach. *Phys. Rev. D* **66** 63006/1-14. Available at arxiv.org as astro-ph/0206501.

HANSEN, F. K., MARINUCCI, D. and VITTORIO, N. (2003). Extended empirical process test for non-Gaussianity in the CMB, with an application to non-Gaussian inflationary models. *Phys. Rev. D* **67** 123004/1-7. Available at arxiv.org as astro-ph/0302202.

JOHNSON, N. L. and KOTZ, S. J. (1972). *Distributions in Statistics*: *Continuous Multivariate Distributions.* Wiley, New York. MR418337

KOMATSU, E. and SPERGEL, D. N. (2001). Acoustic signatures in the primary microwave background bispectrum. *Phys. Rev. D* **63** 63002/1-13. Available at arxiv.org as astro-ph/0005036.

LEONENKO, N. N. (1999). *Limit Theorems for Random Fields with Singular Spectrum.* Kluwer, Dordrecht. MR1687092

LIBOFF, R. L. (1998). *Introductory Quantum Mechanics*, 3rd ed. Addison-Wesley, Reading, MA.

MARINUCCI, D. and PICCIONI, M. (2003). The empirical process on Gaussian spherical harmonics. Working Paper n. 7, Dip. Matematica, Università di Roma "La Sapienza."

MILLER, C. J., NICHOL, R. C., GENOVESE, C. and WASSERMAN, L. (2002). A nonparametric analysis of the cosmic microwave background power spectrum. *Astrophys. J.* **565** L67–L70. Available at arxiv.org as astro-ph/0112049.

Novikov, D., Schmalzing, J. and Mukhanov, V. F. (2000). On non-Gaussianity in the cosmic microwave background. *Astronomy and Astrophysics* **364** 17–25. Available at arxiv.org as astro-ph/0006097.

Peacock, J. A. (1999). *Cosmological Physics*. Cambridge Univ. Press.

Peebles, P. J. E. (1993). *Principles of Physical Cosmology*. Princeton Univ. Press. MR1216520

Phillips, N. G. and Kogut, A. (2001). Statistical power, the bispectrum and the search for non-Gaussianity in the cosmic microwave background anisotropy. *Astrophys. J.* **548** 540–549. Available at arxiv.org as astro-ph/0010333.

Shorack, G. (1972). Convergence of quantile and spacings processes with applications. *Ann. Math. Statist.* **43** 1400–1411. MR359133

Shorack, G. and Wellner, J. (1986). *Empirical Processes with Applications to Statistics*. Wiley, New York. MR838963

van der Vaart, A. and Wellner, J. (1996). *Weak Convergence and Empirical Processes*. Springer, New York. MR1385671

Varshalovich, D. A., Moskalev, A. N. and Khersonskii, V. K. (1988). *Quantum Theory of Angular Momentum*. World Scientific, Singapore. MR1022665

Wasserman, L., Miller, C., Nichol, B., Genovese, C., Jang, W., Connolly, A., Moore, A., Schneider, J. and the PICA Group (2001). Nonparametric inference in astrophysics. Preprint. Available at arxiv.org as astro-ph/0112050.

Winitzki, S. and Wu, J. H. P. (2000). Inter-scale correlations as measures of CMB Gaussianity. Preprint. Available at arxiv.org as astro-ph/0007213.

Wong, E. (1971). *Stochastic Processes in Information and Dynamical Systems*. McGraw-Hill, New York.

Dipartimento di Matematica  
Università di Roma "Tor Vergata"  
via della Ricerca Scientifica 1  
00133 Roma  
Italy  
e-mail: marinucc@mat.uniroma2.it

Dipartimento di Matematica  
Università di Roma "La Sapienza"  
Piazzale Aldo Moro 2  
00185 Roma  
Italy  
e-mail: piccioni@mat.uniroma1.it